\documentclass[submission]{pmsart}[2006/01/13]

\usepackage{amsfonts}

\usepackage{amsmath}

\newcommand{\halmos}{{\mbox{\, \vspace{3mm}}} \hfill
\mbox{$\Box$}}
\newcommand{\bearno}{\begin{eqnarray*}}
\newcommand{\enarno}{\end{eqnarray*}}

\newcommand{\E}{\mathbb{E}}

\newcommand{\R}{\mathbb{R}}
\newcommand{\p}{\mathbb{P}}

\newcommand{\mean}{\mathbb E}
\newcommand{\prob}{\mathbb P}
\newcommand{\de}{\,{\rm d}}

\volume{0}%
\fasc{0}%
\years{0000}%
\pages{000--000}%
\firstpage{1} %
\lastrevision{xx.xx.xxxx}
\title{Two-dimensional ruin probability for subexponential claim size}%                                        %title of the paper
\titlethanks{This work is partially supported by the National Science Centre
under the grants 2013/09/B/HS4/01496 and 2015/17/B/ST1/01102 (2016-2019).
Second and third authors kindly acknowledge partial support
by the project RARE-318984, a Marie Curie IRSES Fellowship
within the 7th European Community Framework Programme.}%                                  %information about support grants etc.
\dedicated{}%                                    %e.g. Dedicated to the memory of Albert Einstein
\ShortTitle{Two-dimensional ruin probability} %                                  %abbreviated title
\ShortAuthors{S. Foss, D. Korshunov, Z. Palmowki and T. Rolski}%                                 %abbreviated authors' names,
\NumberOfAuthors{4}%                              %number of authors of the paper

\received{10.02.2017}%                                     %date of submission

\FirstAuthor{Sergey Foss} %                               %information about first author
\FirstAuthorCity{Edinburgh and Novosibirsk} %
\FirstAuthorAffiliation{Sergey Foss\\Heriot-Watt University, UK\\
Novosibirsk State University, Russia\\
Sobolev Institute of Mathematics, Russia}%
\FirstAuthorAddress{Edinburgh EH14 4AS, Scotland, UK}%
\FirstAuthorEmail{Serguei.Foss@hw.ac.uk}%
\FirstAuthorThanks{}%

\SecondAuthor{Dmitry Korshunov} %                               %information about second author
\SecondAuthorCity{Lancaster} %
\SecondAuthorAffiliation{Dmitry Korshunov\\Lancaster University, UK\\Novosibirsk State University, Russia}%
\SecondAuthorAddress{Lancaster, LA1 4YF, United Kingdom}%
\SecondAuthorEmail{d.korshunov@lancaster.ac.uk}%
\SecondAuthorThanks{}%

\ThirdAuthor{Zbigniew Palmowski}%                                  %first author names and surname. Please enclose names with {} or \mbox
\ThirdAuthorCity{Wroc\l aw}%                              %town being represented by first author
\ThirdAuthorAffiliation{Zbigniew Palmowski\\Faculty of Pure and Applied Mathematics\\
Wroc\l aw University of Science and Technology}%                       %first author affiliation
\ThirdAuthorAddress{Wyb. Wyspia\'nskiego 27, 50-370 Wroc\l aw, Poland}%                           %first author address
\ThirdAuthorEmail{zbigniew.palmowski@gmail.com}                              %first author email
\ThirdAuthorThanks{}                             %information about support grants etc.

\FourthAuthor{Tomasz Rolski}%                                %information about fourth author
\FourthAuthorCity{Wroc{\l}aw}%
\FourthAuthorAffiliation{Tomasz Rolski\\Mathematical Institute\\ University of Wroc{\l}aw}%
\FourthAuthorAddress{pl. Grunwaldzki 2/4, 50-384 Wroc{\l}aw}%
\FourthAuthorEmail{Tomasz.Rolski@math.uni.wroc.pl}%
\FourthAuthorThanks{}%

\FifthAuthor{} %                                %information about fifth author
\FifthAuthorCity{} %
\FifthAuthorAffiliation{}%
\FifthAuthorAddress{} %
\FifthAuthorEmail{}%
\FifthAuthorThanks{}%


\MSCcodes{Primary: 60G51, 60G50; Secondary: 60K25.}

\begin{document}
\maketitle
\begin{abstract}
We analyse the asymptotics of ruin probabilities of two
insurance companies (or two branches of the same company) that
divide between them both claims  and premia in some specified
proportions when the initial reserves of both companies tend to
infinity and generic claim size is subexponential.

\end{abstract}

\section{Introduction}\label{sec:intr}

We consider a model of two-dimensional risk process with renewal input:
claims $\vec{\sigma}_n = (\sigma_{n,1},\sigma_{n,2})$ arrive
in a random input at arrival epochs $\{t_n\}$ with interarrival times
$\{\tau_n\}$.
There are two insurance companies, company $i$  has initial capital
$x_i$ and premium rate $p_i$, and covers claims $\sigma_{n,i}$, $i=1$, $2$.
We assume that the two sequences $\{\vec\sigma_n\}$ and $\{\tau_n\}$
are mutually independent and that each of them consists of i.i.d.
random variables.

Let $N(t)$ be the number of claims by time $t\ge 0$,
$$
N(t) = \max \{n \ : \ t_n \le t\},
$$
then
$$
\vec{S}(t)\ =\ (S_1(t),S_2(t))\ :=\ \sum_{n=1}^{N(t)} \vec{\sigma}_n
$$
is the vector of total claim sizes by time $t$. Further, let
$$
\vec{b}(t)\ =\ (b_1(t),b_2(t))\ :=\ \vec{x} + \vec{p}t
$$
be the sum of initial capitals and of total premiums by time $t$,
here $\vec{x} = (x_1,x_2)$ and $\vec{p}=(p_1,p_2)$.
Introduce ruin probabilities of two types:
for any $T\in (0,\infty]$, these are
\begin{eqnarray*}
\psi_\wedge(\vec{x},T) &=& \prob\{S_1(t) > b_1(t)
\mbox{ or }S_2(t) >b_2(t)\mbox{ for some }t\le T\},\\
\psi_\vee(\vec{x},T) &=& \prob\{S_1(t) > b_1(t)
\mbox{ and }S_2(t) >b_2(t)\mbox{ for some }t\le T\}.
\end{eqnarray*}
Here $\psi_{\wedge}(x_1,x_2, T)$ describes the ruin probability
of at least one insurance company, while $\psi_{\vee}(x_1,x_2,T)$
corresponds to the ruin of the both insurance companies.

There are several papers---see e.g.,
%Avram, Palmowski and Pistorius (2008),
Konstantinides and Li \cite{KL2016}
%Tang (2004)
and list of references therein---where it is assumed
that the pairs $(\sigma_{n,1},\sigma_{n,2})$
are i.i.d. and have a multivariate regularly varying distribution.
We like to consider claim sizes with distributions from a more general subexponential class that
includes Pareto and also log-normal and Weibull distributions.
We are unaware of any reasonable concept of subexponentiality here.
So only two extreme options seem to be doable: either\\
(i) $\sigma_{n,1}$ and $\sigma_{n,2}$ are independent or\\
(ii) they are dependent deterministically, say
$$
\sigma_{n,1} = l(\sigma_n)\quad\mbox{ and }
\quad \sigma_{n,1} = \sigma_n - l (\sigma_n)
$$
where $0\le l(x)\le x$
and $\sigma_n$ is the cumulative claim of the client $n$
that is covered by the two insurance companies together.

There are several papers that study related problems in direction (i)---see
e.g. Li et al. \cite{LLT2007}, Chen et al. \cite{CYN2011},
%Zhang and Wang (2012),
Chen, Wang, Wang \cite{CWW2013_1,CWW2013_2}, Wang et al. \cite{add1}, Jiang et al. \cite{add2},
%Chen, Yuen, Ng (2013),
Hu and Jiang \cite{hu}, Lu and Zhang \cite{LZ2016},
see also references therein.
In particular, in the paper by Lu and Zhang \cite{LZ2016},
a uniform asymptotics over finite intervals have been obtained
under some restrictive assumptions on distributions.

%In second direction, there are, in particular, partial results
%by Badila and Palmowski (2014) and by
%Foss, Rolski and Zachary (2010, conference paper).

Avram et al. \cite{APP, APP2} studied boundary crossing probabilities
of a stochastic process with light-tailed increments.
This study was also motivated
by ruin probabilities of two insurance companies with proportional
claims and the steady state distribution of a
tandem queue with two servers (see \cite{Mandjes}).
Similar considerations were done in Badila et. al. \cite{badila} (for mutidimensional risk process)
and in \cite{hu} (for the two-dimensional risk process with constant interest rate).

The key message of this paper is that in both cases results
on the uniform asymptotics may be obtained for a general
class of {\it strong subexponential distributions},
no further restrictions are needed.
We only consider the second direction using approach developed
in a series of works of Denisov, Foss, Korshunov, Palmowski, and Zachary
(see e.g. \cite{DFK, FZ, FPZ, FKZ, Dima}); similar arguments apply in the first
direction.

To make it simple, we assume $l(x)$ to be linear,
$l(x)= lx$ for some $l\in (0,1)$. Then, after some
transformations we can reformulate the problem as follows. Let us define
\begin{equation}\label{St}
S_t=\sum_{i=1}^{N_t}\sigma_i,
\end{equation}
where $N_t$ is a renewal process with positive i.i.d.
inter-arrival times $\tau_i$, and the claim sizes $\sigma_i$ are
positive i.i.d. random variables that do not depend on $N(t)$ and have a
common distribution function $F(x)$.
%We shall denote by $\lambda$ and
%$\mu$ the reciprocals of the means of $\tau_i$  and $\sigma_i$,
%respectively.
Note also that $N_t=\max\{k\ge 0: t_k\le t\}$ for a random walk constructed from the interarrival times:
$$
t_n=\sum_{k=1}^n\tau_k.
$$

Let the boundaries $b_1$ and $b_2$ be given by
$$
b_1(t) = b_1(t;x_1) = x_1 + p_1t,\quad\quad b_2(t) = b_2(t;x_2) =
x_2 + p_2t.
$$
We assume that $\E \sigma<\infty$, $\E \tau<\infty$ and that
\begin{equation}\label{p1p2}
p_1 > p_2,\quad\quad\quad
p_2>\rho:=\E[\sigma]/\E[\tau]
\end{equation}
for generic $\tau$ and $\sigma$.
%Let $\{\mathcal{F}_t\}_{\{t\ge 0\}}$ be a natural filtration generated by $S_t$ satisfying usual conditions and $\mathcal{T}$ a family of all $\mathcal{F}_t$-stopping times.
In this paper for $T>0$,
%$\varrho\in \mathcal{T}$
we will consider the following boundary crossing probabilities:
\begin{eqnarray*}
\psi_{\wedge}(x_1,x_2,T) &=& \prob\{\tau_\wedge(x_1,x_2) \le T\},\\ \psi_{\vee}(x_1,x_2, T) &=& \prob\{\tau_\vee(x_1,x_2)\le T\},\\
%\psi_{\times}(x_1,x_2) &=& P(\exists t\ge0, u\ge0: S_t > b_1(t),
%S_u
%> b_2(u)),
\end{eqnarray*}
for
\begin{eqnarray}
\tau_\wedge(x_1,x_2)&=&\inf\{t\ge 0: S_t > b_1(t)\wedge
b_2(t)\},\label{tauruin1}\\
\tau_\vee(x_1,x_2)&=&\inf\{t\ge 0: S_t > b_1(t)\vee b_2(t)\},\label{tauruin2}
\end{eqnarray}
where $x\wedge y=\min\{x,y\}$ and $x\vee y = \max\{x,y\}$.

The $\psi_{\wedge}(x_1,x_2,T)$ describes the ruin probability of at
least one insurance company before time $T$ (if $T=\infty$ then
$$
\psi_{\wedge}(x_1,x_2):=\psi_{\wedge}(x_1,x_2,\infty)
$$
is a ruin probability of at least one company).
The $\psi_{\vee}(x_1,x_2,T)$
and
$$\psi_{\vee}(x_1,x_2):=\psi_{\vee}(x_1,x_2,\infty)$$
correspond to the ruin of both insurance companies.
%Finally, $\psi_{\times}(x_1,x_2)$ describes probability that both companies will have ruin.
The first assumption in (\ref{p1p2}) means that the
second company has a smaller premium rate than the first company, and the
second assumption is the condition under which reserves of both
insurance companies tend to infinity. The solutions to the
"degenerate two-di\-men\-si\-o\-nal" ruin problems strongly depend
on the relative position of the vector of premium rates $p=(p_1,
p_2)$ with respect to the proportions vector $(1,1)$. Namely, if
the initial reserves satisfy $x_2 \le x_1$, the two lines do not intersect.
It follows therefore that  $\psi_\wedge$ and $\psi_\vee$
are ruin probabilities of the second and first companies respectively.
In this case the asymptotics follows from
one-dimensional ruin theory---see e.g. Rolski et al. \cite{Rolski}.
Therefore we focus here on the opposite case where
\begin{equation}\label{assump0}
x_1<x_2.
\end{equation}

In this paper we derive the exact first order asymptotics of these
ruin probabilities if $x_1, x_2$ tend to infinity and when the
claims follow a subexponential distribution. We model the claims
by subexponential distributions since many catastrophic events
like earthquakes, storms, terrorist attacks etc are used in their
description. Insurance companies use e.g. the lognormal
distribution (which is subexponential) to model car claims---see Foss et al. \cite{FKZ},
Rolski et al. \cite{Rolski} or Embrechts et al. \cite{Embrechts} for the further
background.

The paper is organized as follows. In next section we present the main results which will be proved in
Section \ref{sec:proofs}.

\section{Main results}\label{main}

In order to state our results we recall basic notation and notions.
Hereinafter we write $f(x,y)\sim g(x,y)$ if $f(x,y)/g(x,y)\to 1$ as $x$, $y\to\infty$
and $f(x)\sim g(x)$ if $f(x)/g(x)\to1$ as $x\to\infty$.
For a distribution $F$, $\overline{F}$ denotes the tail distribution
function given by $\overline{F}(x)=1-F(x)$.

A distribution $F$ on $\R^+$ is called subexponential ($F\in\mathcal{S}$)
if and only if $\overline{F}(x)>0$ for all $x$ and
\begin{equation}\label{eq:5}
\overline{F^{*2}}(x) \ \sim\ 2\overline{F}(x)\quad\text{as }x \to \infty,
\end{equation}
where $F^{*2}$ is the convolution of $F$ with itself.

A distribution $F$ on $\R^+$ is called {\it strong subexponential}
($F\in\mathcal{S}^*$) if and only if $\overline{F}(x)>0$ for all $x$ and
\begin{equation}\label{eq:6}
  \int_0^x\overline{F}(x-y)\overline{F}(y)\,dy
  \sim 2m_{F}\overline{F}(x)\quad\text{as }x \to \infty,
\end{equation}
where
\begin{displaymath}
  m_{F} = \int_0^\infty \overline{F}(x) \,dx
\end{displaymath}
is the expectation of $F$.  It is known that the property $F\in\mathcal{S}^*$
is a tail property of $F$, namely, if $F_1\in\mathcal{S}^*$ and
$\overline{F}_1(x)\sim\overline{F}_2(x)$ as $x\to\infty$, then
$F_2\in\mathcal{S}^*$.  Further, if $F\in\mathcal{S}^*$ then $F\in\mathcal{S}$
and also $F^s\in\mathcal{S}$ where
\begin{displaymath}
  \overline{F^s}(x) = \min \biggl( 1,
    \int_x^{\infty} \overline{F}(t)\, dt
  \biggr)
\end{displaymath}
is the \emph{integrated tail distribution} function
determined by $F$; see \cite{FKZ} for details.

Denote
\begin{equation}\label{ms}
m_i\ =\ p_i \E \tau -\E \sigma,\qquad i=1,2.
\end{equation}
Since $p_1>p_2$, we have $m_1>m_2$.

\begin{theorem}\label{Theorem1}
Assume that the distribution $F$ of generic $\sigma \ge 0$
is strong subexponential. Then, as $x_1$, $x_2\to\infty$,
\begin{equation}\label{min}
\psi_\wedge (x_1,x_2,T)\ \sim\ H_T(x_1,x_2)
\ :=\ \int_0^{\E N_T}
\overline F(\min\{x_1 +t m_1, x_2 +t m_2\})\de t
\end{equation}
and
\begin{equation}\label{max}
\psi_\vee (x_1,x_2, T)\ \sim\ U_T(x_1,x_2)
\ :=\ \int_0^{\E N_T} \overline F(\max\{x_1 +t m_1 , x_2 +t m_2 \})\de t
\end{equation}
holds uniformly for all $T>0$.
\end{theorem}

\begin{corollary}\label{Corollary2}
In conditions of Theorem \ref{Theorem1}, as $x_1$, $x_2$ and $T\to\infty$,
\begin{equation}\label{min1.inf}
\psi_\wedge (x_1,x_2,T)\ \sim\ \int_0^{T/\E\tau}
\overline F(\min\{x_1+t m_1, x_2+t m_2\})\de t,
\end{equation}
and
\begin{equation}\label{max1.inf}
\psi_\vee (x_1,x_2, T)\ \sim\ \int_0^{T/\E\tau}
\overline F(\max\{x_1+t m_1, x_2+t m_2 \})\de t.
\end{equation}
\end{corollary}

\begin{corollary}\label{Corollary1}
Assume that $F$ is strong subexponential and $a<1$ is fixed.
Then we have the following asymptotics, as $x\to\infty$:
\begin{equation}\label{min2}
\psi_\wedge (ax,x)\ \sim\ H(x)\ :=\
\int_0^\infty \overline F(\min\{ax+t m_1, x+t m_2 \})\de t,
\end{equation}
and
\begin{equation}\label{max2}
\psi_\vee (ax,x)\ \sim\ U(x)\ :=\
\int_0^\infty \overline F(\max\{ax +t m_1, x +t m_2\})\de t.
\end{equation}
\end{corollary}

We can also identify the joint distribution of the ruin times \eqref{tauruin1}
and \eqref{tauruin2} and the position at the moment at these ruin times.
To do this we will require additional assumptions.
Let $e(x)$ be a function tending to infinity as $x\to\infty$.
For fixed $y\ge 0$ and $v\ge 0$ we denote:
\begin{equation}\label{HKyv}
H^{y,v}(x)= \int_{ye(x)}^\infty \overline F(\min\{ax+t m_1, x+t m_2\}+ve(x))\de t
\end{equation}
and
\begin{equation}\label{UKyv}
U^{y,v}(x)= \int_{ye(x)}^\infty \overline F(\max\{ax+t m_1, x+t m_2\}+ve(x))\de t.
\end{equation}
Note that $H(x)=H^{0,0}(x)$ and $U(x)=U^{0,0}(x)$.
We will assume that there exists some continuous distribution function $G$ such that
\begin{equation}\label{condition1}
\lim_{x\to\infty}\frac{H^{y,v}(x)}{H(x)}=\overline{G}(m_1y+v)
\end{equation}
and that
\begin{equation}\label{condition2}
\lim_{x\to\infty}\frac{U^{y,v}(x)}{U(x)}=\overline{G}(m_2y+v).
\end{equation}

\begin{theorem}\label{Theorem2}
Assume that there exist a function $e(x) \uparrow\infty $ and a
continuous probability distribution $G$ on the positive half-line $(0,\infty)$ such that
\eqref{condition1} and \eqref{condition2} are satisfied.
Then, for $a<1$, we have:
\begin{equation}\label{min2}
\lim_{x\to\infty}\prob\left\{ \frac{\tau_\wedge(ax,x)}{e(x)}\ge y, \quad \frac{S_{\tau_\wedge(ax,x)}-b_\wedge(t)}{e(x)}\ge v\Bigg|\tau_\wedge(ax,x)<\infty\right\}=\overline{G}(m_1y+v)
\end{equation}
and
\begin{equation}\label{max2}
\lim_{x\to\infty}\prob\left\{ \frac{\tau_\vee(ax,x)}{e(x)}\ge y, \quad  \frac{S_{\tau_\vee(ax,x)}-b_\vee(t)}{e(x)}\ge v\Bigg|\tau_\vee(ax,x)<\infty\right\}=\overline{G}(m_2y+v).
\end{equation}
\end{theorem}

For similar statements for a linear boundary,
see Asmussen and Kl\"uppelberg \cite{AK1996} and
Foss et al. \cite{FKZ} in the case of i.i.d. jumps
and Asmussen and Foss \cite{AF2014} in the case of Markov modulation.

\begin{remark}\rm
To provide examples where
assumptions \eqref{condition1} and \eqref{condition2} are satisfied, we
choose strong subexponential distribution $F$ for which $F^s$ is self-neglecting. That is,
let us assume that there exist a function $e(x) \uparrow\infty $ and a
continuous probability distribution $G$ on the positive half-line $(0,\infty)$ such that, for
any $y>0$,
\begin{equation}\label{defF*}
\frac{\overline{F}^s(x+ye(x))}{\overline{F}^s(x)} \to\ \overline{G}(y)\quad \text{as }x\to\infty.
\end{equation}

From \cite[Thm. 3.4.5, p. 158]{Embrechts} it follows that $F^s$ is either in the domain
of attraction of Frechet distribution
or in the the domain of attraction attraction of Gumbel distribution.
In the first case
\begin{equation}\label{Frechet}
F\in {\rm RV}(\alpha+1),\quad F^s\in {\rm RV}(\alpha),\quad \overline{G}(y)=(1+\alpha^{-1}y)^{-\alpha},\quad e(x)=\alpha^{-1}x,
\end{equation}
with $\alpha >0$ (then $\E \sigma <\infty$) for the family RV$(\alpha+1)$ of regularly varying distributions.
In the second case,
\begin{equation}\label{Gumbel}
\overline{G}(y)=e^{-y},\qquad e(x) =\frac{m_F\overline{F}^s(x)}{\overline{F}(x)}.
\end{equation}

To check \eqref{condition1} it is convenient to use the following representation
\begin{eqnarray}
H^{y,v}(x)&=&\frac{1}{m_1} \overline{F}^s(ax+(m_1y+v)e(x))\nonumber\\
&&+\left(\frac{1}{m_2}-\frac{1}{m_1}\right)
\overline{F}^s(ax+m_1\widehat T+ve(x))\nonumber\\&&
+\frac{1}{m_2}\overline{F}^s(x+m_2\widehat T+ve(x))\label{representation}\end{eqnarray}
for the moment when two lines $ax+m_1t$ and $x+m_2t$ cut each other:
\[
\widehat T\ =\ \frac{(1-a)x}{m_1-m_2}.
\]
Recall also that $H(x)=H^{0,0}(x)$. Now one can check that if \eqref{Frechet} holds true then
indeed \eqref{condition1} is satisfied. Similarly, one can prove that in this case
assumption \eqref{condition1} is also satisfied.

The case \eqref{Gumbel} is much more advanced and it should be analysed
case by case using for example \cite[Example 3.3.35, p. 149]{Embrechts}.
\end{remark}

In our proof of Theorem \ref{Theorem1}, we use the following simple
result that may be of use in other settings.

\begin{lemma}\label{l:int}
Let $z(t)\ge 0$ be an increasing function, $F$ a distribution
and $N_t$ a renewal process. Then
\begin{eqnarray}\label{Eint.upper}
\E \int_0^{N_T} \overline{F}(x+z(t)) dt
&\le& \int_0^{\E N_T} \overline{F}(x+z(t)) dt
\quad\mbox{for all }x\mbox{ and }T>0
\end{eqnarray}
and
\begin{eqnarray}\label{Eint.asy.Tinf}
\E \int_0^{N_T} \overline{F}(x+z(t)) dt
&\sim& \int_0^{\E N_T} \overline{F}(x+z(t)) dt
\quad\mbox{as }x,\ T\to\infty.
\end{eqnarray}
If, in addition, the distribution $F$ is a long-tailed, then
\begin{eqnarray}
\label{Eint.asy}&\\
\E \int_0^{N_T} \overline{F}(x+z(t)) dt
&\sim& \int_0^{\E N_T} \overline{F}(x+z(t)) dt
\quad\mbox{as }x\to\infty\mbox{ uniformly for all }T>0.\nonumber
\end{eqnarray}
\end{lemma}

Notice that, in general, the last result is not applicable to an upper
limit $N_\varrho$ where $\varrho$ is a stopping time like in \eqref{rho1}
and \eqref{rho2} below.

\section{Proofs}\label{sec:proofs}

\subsection{Proof of Theorem \ref{Theorem1}}
As discussed above, it is sufficient to consider the case $x_2>x_1$.
For any $x>0$ and $\widehat T>0$, let $\widehat x=x+p_1\widehat T$ and
\begin{eqnarray*}
\widehat b(t) &:=& \left\{
\begin{array}{ll}
x+p_1t &\mbox{ for }t\le\widehat T\\
\widehat x+p_2(t-\widehat T) &\mbox{ for }t>\widehat T,
\end{array}
\right.
\end{eqnarray*}
which is a continuous piece-wise linear function.
Then both assertions of the theorem, \eqref{min} and \eqref{max},
will follow if the following tail asymptotics is proven
for all $p_1>0$ and $p_2>0$, without assumption $p_2<p_1$:
\begin{eqnarray}\label{psi.gen}
\psi(x,\widehat T,T) &:=& \prob\{S_t>\widehat b(t)\mbox{ for some }t\le T\}\nonumber\\
&\sim& \int_0^{\E N_T}\overline F(\widehat x(t))dt
\quad\mbox{as }x\to\infty\mbox{ uniformly for all }\widehat T\mbox{ and }T,
\end{eqnarray}
where
\begin{eqnarray*}
\widehat x(t) &:=& \left\{
\begin{array}{ll}
x+m_1t &\mbox{ for }t\le\widehat T/\E\tau\\
x+m_1\widehat T/\E\tau+m_2(t-\widehat T/\E\tau) &\mbox{ for }t>\widehat T/\E\tau.
\end{array}
\right.
\end{eqnarray*}

Since $F\in{\mathcal S}^*$, $F$ is particularly subexponential
and long-tailed.
For any fixed $T$, the random variable $N_T$ has a light tail,
that is, it possesses a finite exponential moment.
In addition, $N_T$ is independent of $\sigma$'s.
Hence, by Theorem 3.37 in \cite{FKZ}, for every fixed $T$,
\begin{eqnarray}
&&\label{as.for.eta}\\\nonumber
\p\Bigl\{\sum_{i:t_i\le T}\sigma_i>x\Bigr\}
&=& \sum_{j=0}^\infty \p\{N_T=j\}\overline{F^{*j}}(x)
\ \sim\  \E N_T\overline F(x)\quad\mbox{ as }x\to\infty,
\end{eqnarray}
and this equivalence holds uniformly on any $T$-compact set.
This observation together with the following lower and upper bounds
\begin{eqnarray*}
\p\Bigl\{\sum_{i:t_i\le T}\sigma_i>x+\max\{p_1,p_2\}T\Bigr\}
&\le& \psi(x,\widehat T,T)\ \le\
\p\Bigl\{\sum_{i:t_i\le T}\sigma_i>x\Bigr\}
\end{eqnarray*}
and long-tailedness of $F$ yields that,
uniformly on any $T$-compact set,
\begin{eqnarray*}
\psi(x,\widehat T,T) &\sim& \E N_T\overline F(x)
\quad\mbox{ as }x\to\infty.
\end{eqnarray*}
Taking into account that
\begin{eqnarray*}
N_T\overline F(x+\max\{m_1,m_2\}T)\ \le\
\int_0^{N_T}\overline F(\widehat x(t))dt
\ \le\ N_T\overline F(x)
\end{eqnarray*}
and that $F$ is long-tailed, we prove that \eqref{psi.gen}
holds as $x\to\infty$ uniformly on $T$-compact sets.

Therefore, it remains to prove \eqref{psi.gen} for the case $T\to\infty$.
If $T\le\widehat T$, then
\begin{eqnarray*}
\psi(x,\widehat T,T) &=& \prob\Bigl\{\sum_{j=1}^n\sigma_j>x+p_1t_n
\mbox{ for some }n\le N_T\Bigr\},
\end{eqnarray*}
so, as proven in \cite[Theorem 3]{Dima} for the supremum of
a compound renewal process with negative drift
$\E\sigma/\E\tau-p_1=-m_1/\E\tau$,
as $x\to\infty$ uniformly for $T\le\widehat T$,
\begin{eqnarray}\label{psi.wedge.1}
\psi(x,\widehat T,T) &\sim& \frac{1}{m_1}
\int_0^{m_1\E N_T}\overline F(x+t)dt\\
&=& \int_0^{\E N_T}\overline F(x+m_1 t)dt
\ \sim\ \int_0^{\E N_T}\overline F(\widehat x(t))dt,
\end{eqnarray}
where the last equivalence follows because $\E N_T\sim T/\E\tau$
as $T\to\infty$ and, for any $\varepsilon>0$,
\begin{eqnarray}\label{int.N.x}
\lefteqn{(1-\varepsilon)
\int_0^{T/\E\tau}\overline F(x+m_1t)dt
\le \int_0^{(1-\varepsilon)T/\E\tau}\overline F(x+m_1t)dt}\nonumber\\
&&\le \int_0^{(1+\varepsilon)T/\E\tau}\overline F(x+m_1t)dt
\ \le\ (1+\varepsilon)
\int_0^{T/\E\tau}\overline F(x+m_1t)dt,
\end{eqnarray}
and similar bounds for the integral of $\overline F(\widehat x(t))$.

Let us now consider the case $T\to\infty$ and $T>\widehat T$.
If $\widehat T$ is bounded, then we make use of the inequalities
\begin{eqnarray*}
\lefteqn{\prob\{S_t>x+p_1\widehat T+p_2t\mbox{ for some }t\le T\}
\le \psi(x,\widehat T,T)}\\
&&\le \prob\{S_t>x-p_1\widehat T+p_2t\mbox{ for some }t\le T\},
\end{eqnarray*}
which reduce---due to long-tailedness of $F$---the problem to the case
$\widehat T=T$, a particular case of considered above.

The case where $T\to\infty$, $\widehat T\to\infty$ and $T>\widehat T$,
but $T-\widehat T$ is bounded, is very similar.
Indeed, in this case it is enough to notice that
\begin{eqnarray*}
\lefteqn{\prob\{S_t>x+p_2(T-\widehat T)+p_1t\mbox{ for some }t\le T\}
\le \psi(x,\widehat T,T)}\\
&&\le \prob\{S_t>x-p_2(T-\widehat T)+p_1t\mbox{ for some }t\le T\}.
\end{eqnarray*}

Let us now consider the last remaining case where $T\to\infty$,
$\widehat T\to\infty$, and $T-\widehat T\to\infty$.
We demonstrate two approaches, where the first one is based on the uniform
equivalences in the case of linear functions obtained in \cite{Dima}
and the second on the discrete time results from \cite{FPZ}.

\subsection{Proof based on \cite{Dima}}
Since $F$ is particularly long-tailed, there exists an
increasing function $h(x)\to\infty$ such that $h(x)=o(x)$ and
\begin{eqnarray}\label{hx}
\overline F(x+h(x)) &\sim& \overline F(x)\quad\mbox{as }x\to\infty,
\end{eqnarray}
see \cite[Lemma 2.19]{FKZ}. For $T>\widehat T$,
\begin{eqnarray}
\nonumber\lefteqn{\psi(x,\widehat T,T) \le
\prob\Bigl\{\sum_{j=1}^n\sigma_j>x+p_1t_n-h(x)
\mbox{ for some }n\le N_{\widehat T}\Bigr\}}\\
&&
%\hspace{5mm}
+\prob\Bigl\{\sum_{j=1}^n\sigma_j\le x+p_1t_n-h(x)
\mbox{ for all }n\le N_{\widehat T},\nonumber\\&&\qquad\qquad
\sum_{j=1}^n\sigma_j>\widehat x(t_n)
\mbox{ for some }n\in(N_{\widehat T},N_T]\Bigr\}\nonumber\\
&&= P_1+P_2.\label{max.u}
\end{eqnarray}
Again by \cite[Theorem 3]{Dima},
\begin{eqnarray}\label{max.u.1}
\quad P_1 &\sim&
\int_0^{\E N_{\widehat T}}\overline F(x-h(x)+m_1t)dt
\ \sim\ \int_0^{\E N_{\widehat T}}\overline F(x+m_1t)dt
\quad\mbox{as }x\to\infty,
\end{eqnarray}
as follows from \eqref{hx}.

The second probability on the right hand side of \eqref{max.u} equals
\begin{eqnarray}
P_2 &=& \prob\Bigl\{S_{\widehat T}\le x+p_1\widehat T-h(x),\nonumber\\
\label{P.2.vee} && S_{\widehat T}-p_1\widehat T
+\sum_{j=N_{\widehat T}+1}^n\sigma_j-p_2(t_n-\widehat T)>x
\mbox{ for some }n\in[N_{\widehat T}+1,N_T]\Bigr\}.
\end{eqnarray}
We need to exclude dependence caused by the interval
$(\widehat T,t_{N_{\widehat T}+1}]$. We do it in the following way:
\begin{eqnarray*}
&&\sup_{n\in[N_{\widehat T}+1,N_T]}
\biggl(\sum_{j=N_{\widehat T}+1}^n\sigma_j-p_2(t_n-\widehat T)\biggr)
\le \sigma_{N_{\widehat T}+1}\\
&&\qquad +\sup_{n\in[N_{\widehat T}+2,N_T]}
\biggl(\sum_{j=N_{\widehat T}+2}^n\sigma_j-p_2(t_n-t_{N_{\widehat T}+1})\biggr),
\end{eqnarray*}
where the two random variables on the right are independent because
$t_{N_{\widehat T}+1}$ is a stopping time.
Further, extending the interval $(t_{N_{\widehat T}+1},T]$ to the interval
$(t_{N_{\widehat T}+1},t_{N_{\widehat T}+1}+T-\widehat T]$
of length $T-\widehat T$ we conclude that
\begin{eqnarray*}
\sup_{n\in[N_{\widehat T}+2,N_T]}
\biggl(\sum_{j=N_{\widehat T}+2}^n\sigma_j-p_2(t_n-t_{N_{\widehat T}+1})\biggr),
\end{eqnarray*}
is stochastically not greater than
\begin{eqnarray*}
\zeta_{T-\widehat T} &:=& \sup_{t\le T-\widehat T} (S_t-p_2t),
\end{eqnarray*}
whose tail is equivalent to, again by \cite[Theorem 3]{Dima},
\begin{eqnarray*}
\prob\{\zeta>z\} &\sim& \int_0^{\E N_{T-\widehat T}}\overline F(z+m_2y)dy
\quad\mbox{as }z\to\infty.
\end{eqnarray*}
Fix any $\varepsilon>0$ such that $m_1=p_1\E\tau-\E\sigma>\varepsilon$.
Then it follows from \eqref{P.2.vee} that,
with $\widetilde\zeta_{T-\widehat T}$ being
an independent copy of $\zeta_{T-\widehat T}$,
\begin{eqnarray}\label{est.p2vee}
P_2 &\le& \int_{-\infty}^{x-h(x)}
\prob\{S_{\widehat T}-p_1\widehat T\in dy\}
\prob\{\sigma_{N_{\widehat T}+1}+\widetilde\zeta_{T-\widehat T}>x-y\}\nonumber\\
&=& \int_{-\infty}^{x-h(x)}
\prob\{S_{\widehat T}-p_1\widehat T+(m_1-\varepsilon)\widehat T/\E\tau
\in dy+(m_1-\varepsilon)\widehat T/\E\tau\}\nonumber\\
&&\qquad\prob\{\sigma_{N_{\widehat T}+1}+\widetilde\zeta_{T-\widehat T}>x-y\}\nonumber\\
&=& \int_{-\infty}^{x+(m_1-\varepsilon)\widehat T/\E\tau-h(x)}
\prob\{S_{\widehat T}-(\E\sigma+\varepsilon)\widehat T/\E\tau\in dy\}\nonumber\\
&&\qquad\prob\{\sigma_{N_{\widehat T}+1}+\widetilde\zeta_{T-\widehat T}
>x+(m_1-\varepsilon)\widehat T/\E\tau-y\}.
\end{eqnarray}
We have
\begin{eqnarray*}
\prob\{\sigma_{N_{\widehat T}+1}+\widetilde\zeta_{T-\widehat T}>u\}
&\sim& \prob\{\widetilde\zeta_{T-\widehat T}>u\}+\overline F(u)
\quad\mbox{as }u\to\infty,
\end{eqnarray*}
by \cite[Corollary 3.16]{FKZ}. Since the process
$S_t-(\E\sigma+\varepsilon)t/\E\tau$ has negative drift
$-\varepsilon/\E\tau$ and its value at time $\widehat T$
does not exceed the value of its maximum at the same time,
we also have
\begin{eqnarray*}
\prob\{S_{\widehat T}-(\E\sigma+\varepsilon)\widehat T/\E\tau>u\}
&\le& c_1(\varepsilon)\prob\{\zeta_{\widehat T}>u\},
\end{eqnarray*}
where
\begin{eqnarray*}
\zeta_{\widehat T} &:=& \sup_{t\le \widehat T} (S_t-p_1t).
\end{eqnarray*}
Integrating \eqref{est.p2vee} by parts, applying the last upper bound
and then integrating by parts back we deduce that
\begin{eqnarray*}
\lefteqn{\int_{h(x_1)}^{x+(m_1-\varepsilon)\widehat T/\E\tau-h(x)}
\prob\{S_{\widehat T}-(\E\sigma+\varepsilon)\widehat T/\E\tau\in dy\}}\\&&\qquad
\prob\{\sigma_{N_{\widehat T}+1}+\widetilde\zeta_{T-\widehat T}
>x+(m_1-\varepsilon)\widehat T/\E\tau-y\}\\
&&=\ o\Bigl(\prob\{\zeta_{T-\widehat T}>x+(m_1-\varepsilon)\widehat T/\E\tau\}
+\prob\{\zeta_{\widehat T}>x+(m_1-\varepsilon)\widehat T/\E\tau\}\Bigr),
\end{eqnarray*}
by \cite[Theorem 3.28]{FKZ}, because $F\in\mathcal S^*$. Hence,
\begin{eqnarray*}
P_2 &\le& \int_{-\infty}^{h(x)}
\prob\{S_{\widehat T}-(\E\sigma+\varepsilon)\widehat T/\E\tau\in dy\}\\&&\qquad
\prob\{\sigma_{N_{\widehat T}+1}+\widetilde\zeta_{T-\widehat T}
>x+(m_1-\varepsilon)\widehat T/\E\tau-y\}\\
&&+o\Bigl(\prob\{\zeta_{T-\widehat T}>x+(m_1-\varepsilon)\widehat T/\E\tau\}
+\prob\{\zeta_{\widehat T}>x+(m_1-\varepsilon)\widehat T/\E\tau\}\Bigr)\\
&\le& \prob\{\sigma_{N_{\widehat T}+1}+\widetilde\zeta_{T-\widehat T}
>x+(m_1-\varepsilon)\widehat T/\E\tau-h(x)\}
+o(...)\\
&\sim& \prob\{\widetilde\zeta_{T-\widehat T}>x+(m_1-\varepsilon)\widehat T/\E\tau\}
+o(\prob\{\zeta_{\widehat T}>x+(m_1-\varepsilon)\widehat T/\E\tau\})\\\\
&\sim& \int_0^{\E N_T-\E N_{\widehat T}}
\overline F(x+(m_1-\varepsilon)\widehat T/\E\tau+m_2t)dt\\&&\qquad
+o(\prob\{\zeta_{\widehat T}>x+(m_1-\varepsilon)\widehat T/\E\tau\})\\
&\le& \int_{(m_1-\varepsilon)\widehat T/\E\tau}^{\E N_T}
\overline F(x+m_2t)dt+o(\prob\{\zeta_{\widehat T}>x+(m_2-\varepsilon)\widehat T/\E\tau\}).
\end{eqnarray*}
for all sufficiently large $\widehat T$.
Being substituted together with \eqref{max.u.1}
into \eqref{max.u} it implies that
\begin{eqnarray*}
\psi(x,\widehat T,T) &\le&
(1+o(1))\biggl(\int_0^{\E N_{\widehat T}}\overline F(x+m_1t)dt
+\int_{(m_1-\varepsilon)\widehat T/\E\tau}^{\E N_T}
\overline F(x+m_2t)dt\biggr).
\end{eqnarray*}
Since $\E N_{\widehat T}\sim \widehat T/\E\tau$ as $\widehat T\to\infty$
and due to \eqref{int.N.x} with $T=\widehat T$, we get
\begin{eqnarray*}
\psi(x,\widehat T,T) &\le&
(1+o(1))\biggl(\int_0^{\widehat T/\E\tau}\overline F(x+m_1t)dt
+\int_{\widehat T/\E\tau}^{\E N_T}\overline F(x+m_2t)dt\biggr)
\end{eqnarray*}
as  $x,\ T,\ \widehat T,\ T-\widehat T\to\infty.$
So,
\begin{eqnarray}\label{min.uu}
\qquad\psi(x,\widehat T,T) &\le&
(1+o(1))\int_0^{\E N_T}\overline F(\widehat x(t))dt
\quad\mbox{as } x,\ T,\ \widehat T,\ T-\widehat T\to\infty.
\end{eqnarray}

For the lower bound, let us fix an $\varepsilon>0$
and follow the standard arguments based on the strong law of
large numbers and the single big jump principle.
By the strong law of large numbers, there exists an $A$ such that
\begin{eqnarray}\label{SLLN}
\prob\{|t_n-n\E \tau|<n\varepsilon+A \mbox{ for all }n\ge 1\}
&\ge& 1-\varepsilon.
\end{eqnarray}
On the event introduced in \eqref{SLLN}, if
$n\le \bigl[\frac{\widehat T-A}{\E \tau+\varepsilon}\bigr]=:\widehat n$,
then $t_n\le\widehat T$ and hence $n\le N_{\widehat T}$. Further,
on the event introduced in \eqref{SLLN}, if
$n_1:=\bigl[\frac{\widehat T-A}{\E \tau-\varepsilon}\bigr]
< n\le \bigl[\frac{T-A}{\E \tau+\varepsilon}\bigr]=:n_2$,
then $\widehat T<t_n\le T$ and hence $N_{\widehat T}<n\le N_T$.

Then, since $\bar\sigma$'s do not depend on the renewal process $N_t$,
we obtain the inequality
\begin{eqnarray}\label{psi.lower.pre}
\lefteqn{\psi(x,\widehat T,T) \ge
(1-\varepsilon)\prob \Bigl\{\sum_{i=1}^n \bar\sigma_i>
\widehat b(n(\E \tau+\varepsilon)+A)-n\E\sigma}\\&&
\qquad \mbox{ for some }n\le \widehat n\mbox{ or }n\in(n_1,n_2]\Bigr\}.\nonumber
\end{eqnarray}
Then the standard arguments based on the strong law of
numbers---now for $\sigma$'s---and the single big jump principle finally imply the lower bound
\begin{eqnarray*}
\psi(x,\widehat T,T) &\ge&
(1+o(1))\int_0^{\E N_T}\overline F(\widehat x(t))dt
\quad\mbox{as } x,\ T,\ \widehat T,\ T-\widehat T\to\infty,
\end{eqnarray*}
which together with \eqref{min.uu} concludes the proof.
\halmos

\subsection{Proof based on \cite{FPZ}}
Let us now give an alternative proof of the asymptotic behaviour of
$\psi(x,\widehat T,T)$ in the case $T\to\infty$
based on the results for discrete time from \cite{FPZ}.

Let $S_n=\sum_{k=1}^n \sigma_k$, ${\bar\sigma}_n= \sigma_n-\E\sigma$,
and ${\bar S}_n= \sum_{k=1}^n \bar{\sigma}_k \equiv S_n-n\E\sigma$
be a centered random walk.
From either Corollary 3 of \cite{FZ} or Theorem 1 of \cite{DFK},
it follows that, for any random variable $\gamma$ having a light-tailed
distribution (that is, $\E \exp (\delta\gamma )<\infty$ for some $\delta>0$),
we have, for any real $C$,
\begin{equation}\label{DFKT}
\prob\{{\bar S}_n > x+nC\mbox{ for some }n\le\gamma\}
\ \sim\ \E\gamma \overline F(x).
\end{equation}

Take $\varepsilon>0$ sufficiently small, such that
$p_1(\E\tau-\varepsilon)-\E\sigma>0$ and
$p_2(\E\tau-\varepsilon)-\E\sigma>0$. We take
$$
\gamma=\min\{n\ge 1:\  t_k\ge (\E\tau-\varepsilon)k
\mbox{ for all }k\ge n\},
$$
which is finite almost surely due to the strong law of large numbers.
Moreover, this random variable $\gamma$ is light-tailed because,
for $\delta>0$,
\begin{eqnarray*}
\prob\{\gamma>n\} &=&
\prob\{t_k < (\E\tau-\varepsilon)k \mbox{ for some }k\ge n\}\\
&=& \prob\{e^{\delta((\E\tau-\varepsilon)k-t_k)}>1 \mbox{ for some }k\ge n\}\\
&\le& (\E e^{\delta(\E\tau-\varepsilon-\tau)})^n,
\end{eqnarray*}
by Doob's inequality provided that
$\E e^{\delta(\E\tau-\varepsilon-\tau)}<1$; it holds
for all sufficiently small $\delta>0$, since the random variable
$\E\tau-\varepsilon-\tau$ is bounded above and has negative mean.

Then
\begin{eqnarray*}
&&\prob\{\bar{S}_n>\widehat b(t_n)-n\E\sigma
\mbox{ for some }n\le\gamma\}\\
&&\qquad\le \prob\{\bar{S}_n>x-n\E\sigma\mbox{ for some }n\le\gamma\}
\ \sim\ \E\gamma \overline{F}(x),
\end{eqnarray*}
by \eqref{DFKT}.
Therefore, for all $T>0$,
\begin{eqnarray*}
\psi(x,\widehat T,T) &\le&
\prob\{\bar{S}_n>\widehat b(t_n)-n\E\sigma\mbox{ for some }n\le\gamma\}\\&&+\prob\{\bar{S}_n >
\widehat b(n(\E\tau-\varepsilon))-n\E\sigma
\mbox{ for some }n\in(\gamma,N_T]\}\\
&\le& O(\overline F(x))
+\prob\{\bar{S}_n > \widehat b(n(\E\tau-\varepsilon))-n\E\sigma
\mbox{ for some }n\le N_T\}
\end{eqnarray*}
as $x\to\infty$ uniformly for all $T>0$.

Now we recall the following result where the class $\Gamma$
is a class of all counting random variables $\gamma$ such that,
for all $n$, the event $\{\gamma\le n\}$
does not depend on $\{\sigma_k\}_{k>n}$.
For $c>0$, let $\mathcal{G}_c$ be the class of functions $g$
such that $g(n+1)\ge g(n)+c$ for all $n=1,2,\ldots$.

\begin{theorem}\label{FPZT} (see \cite[Theorem 2(ii)]{FPZ}).
Assume ${\bar S}_n = \sum_1^n \bar{\sigma}_n$
is a centered random walk where the common distribution
of $\bar{\sigma}_n$ belongs to the class $\mathcal{S}^*$.
Then, for any $c>0$, uniformly for all $g\in \mathcal{G}_c$
and for all random times $\gamma\in\Gamma$, we have
\begin{eqnarray}\label{uni}
\prob\{\max_{n\le \gamma} ({\bar S}_n - g(n))>x\} &\sim&
\sum_{n\ge 1} \prob (\gamma \ge n) \prob (\bar{\sigma}_1>x+g(n))\nonumber\\
&\sim& \sum_{n\ge 1} \prob (\gamma \ge n) \overline{F}(x+g(n))
\quad\mbox{as }x\to\infty.
\end{eqnarray}
\end{theorem}

Consider function $g$ of the form
\begin{eqnarray*}
g(n) &=& \widehat b(n(\E\tau-\varepsilon))-n\E\sigma-x,
\end{eqnarray*}
which is in the class $\mathcal G_{p_1\wedge p_2}$.
By Theorem \ref{FPZT}, we get that
\begin{eqnarray*}
\lefteqn{\prob\{\bar{S}_n > \widehat b(n(\E\tau-\varepsilon))-n\E\sigma
\mbox{ for some }n\le N_T\}}\\
&&\qquad \sim \sum_n \prob\{N_T\ge n\} \overline{F}(x+g(n)),
\end{eqnarray*}
as $x\to\infty$ and $T\to\infty$.
Further, since $\varepsilon >0$ is arbitrary,
we may follow the proof of the upper bound for $\psi$ and
let $\varepsilon\to 0$ to obtain an upper bound of the form
$$
(1+o(1)) \sum_n \prob (N_T\ge n)\overline{F}(\widehat x(n)).
$$
The last sum may be rewritten as
\begin{eqnarray*}
\E \sum_n {\mathbb I}(N_T\ge n) \overline{F}(\widehat x(n))
&=& \E \sum_{n=1}^{N_T} \overline{F}(\widehat x(n)) \\
&\sim& \E \int_0^{N_T} \overline{F}(\widehat x(t)) dt
\quad\mbox{as }x\to\infty.
\end{eqnarray*}
Here the last equivalence follows from the long-tailedness of $F$.
Finally, we may use Lemma \ref{l:int} to conclude
with the upper bound \eqref{min.uu}.

Theorem \ref{FPZT} is also applicable to \eqref{psi.lower.pre},
so hence the correct lower bound for $\psi(x,\widehat T, T)$ follows too.
\halmos

\subsection{Proof of Lemma \ref{l:int}}
For any fixed $x$, the function
\begin{eqnarray*}
f(y) &:=& \int_0^y \overline{F}(x+z(t)) dt
\end{eqnarray*}
is a concave function in $y$ because $f'(y)=\overline{F}(x+z(y))$
is decreasing in $y$. Then the upper bound \eqref{Eint.upper}
follows by Jensen's inequality for concave functions.

Concerning \eqref{Eint.asy.Tinf}, notice that,
for any fixed $\varepsilon>0$,
\begin{eqnarray*}
\E \int_0^{N_T} \overline{F}(x+z(t)) dt
&\ge& \prob\{N_T>(1-\varepsilon)\E N_T\}
\int_0^{(1-\varepsilon)\E N_T} \overline{F}(x+z(t)) dt,
\end{eqnarray*}
where $\prob\{N_T>(1-\varepsilon)\E N_T\}\to 1$ as $T\to\infty$
by the law of large numbers for the renewal process and
\begin{eqnarray*}
\int_0^{(1-\varepsilon)\E N_T} \overline{F}(x+z(t)) dt
&\ge& (1-\varepsilon)\int_0^{\E N_T} \overline{F}(x+z((1-\varepsilon)t))dt\\
&\ge& (1-\varepsilon)\int_0^{\E N_T} \overline{F}(x+z(t)) dt,
\end{eqnarray*}
which implies lower bound
\begin{eqnarray*}
\E \int_0^{N_T} \overline{F}(x+z(t)) dt
&\ge& (1+o(1))\int_0^{\E N_T} \overline{F}(x+z(t)) dt
\quad\mbox{as }x,\ T\to\infty,
\end{eqnarray*}
which together with upper bound \eqref{Eint.upper}
justifies \eqref{Eint.asy.Tinf}.

Finally, for any fixed $T$ and $A$, since $z(t)\ge 0$ increases,
\begin{eqnarray*}
\E \int_0^{N_T} \overline{F}(x+z(t)) dt
&\ge& \E\biggl\{\int_0^{N_T} \overline{F}(x+z(t)) dt;\ N_T\le A\biggr\}\\
&\ge& \overline{F}(x+z(A))\E\{N_T;\ N_T\le A\}\\
&\sim& \overline{F}(x)\E\{N_T;\ N_T\le A\}\quad\mbox{as }x\to\infty,
\end{eqnarray*}
provided $F$ is long-tailed in which case also
\begin{eqnarray*}
\int_0^{\E N_T} \overline{F}(x+z(t)) dt &\sim&
\overline{F}(x)\E N_T \quad\mbox{as }x\to\infty,
\end{eqnarray*}
so hence the asymptotics \eqref{Eint.asy}
follows uniformly on $T$-compact sets.
Together with \eqref{Eint.asy.Tinf} it completes the proof.
\halmos

\subsection{Proof of Theorem \ref{Theorem2}}
The proof of this result is very similar to the proof of \eqref{max} of Theorem \ref{Theorem1}.
For example to prove \eqref{min2} it suffices to observe that
\begin{eqnarray*}
\lefteqn{\prob\left\{ \frac{\tau_\wedge(ax,x)}{e(x)}\ge y, \
\frac{S_{\tau_\wedge(ax,x)}-b_\wedge(t)}{e(x)}\ge v;
\tau_\wedge(ax,x)<\infty\right\}}\\
&=&\prob\left\{\bar S_n > \min\{ax +p_1t_n-n\E \sigma,
x+p_2 t_n - n\E \sigma\}+ve(x)\right.\\&&\qquad\qquad \left.\mbox{ for some }n\ge ye(x)\right\}.
\end{eqnarray*}
\halmos

\subsection{Generalisation to $n$ insurance companies}
It is more or less clear that a very similar proof works for \eqref{psi.gen}
if $\widehat x(t)$ is a continuous piece-wise linear increasing function.
Then it allows to go beyond two insurance companies
and consider a model with an arbitrary number of them.
Clearly in high dimension one has to overcome extra
combinatorial problems that seem to be doable.

\subsection{Open question}
We strongly believe that an analogue of Theorem \ref{Theorem1}
holds for stopping times, say $\varrho$ instead of $T$.
We expect the following to be correct:
\begin{equation}\label{rho1}
\psi_\wedge (x_1,x_2,\varrho)\ \sim\
\E\int_0^{N_\varrho}\overline F(\min\{x_1 +t m_1, x_2 +t m_2\})\de t
\end{equation}
and
\begin{equation}\label{rho2}
\psi_\vee (x_1,x_2, \varrho)\ \sim\
\E\int_0^{N_\varrho} \overline F(\max\{x_1 +t m_1 , x_2 +t m_2 \})\de t
\end{equation}
hold as $x\to\infty$ uniformly for all stopping times $\varrho$
with respect to the filtration generated by the renewal process $N_t$.

{\bf Acknowledgements.} We thank Serban Badila for discussions
and preliminary results included in his thesis which motivated us
to work on the subject and to come to much shorter proofs
of more general results given in Theorems \ref{Theorem1} and \ref{Theorem2}.


\begin{thebibliography}{10}
%\bibliographystyle{plain}

\bibitem{AF2014}
Asmussen, S. and Foss, S. (2014)
On exceedance times for some processes with dependent increments.
{\it J. Appl. Probab.} {\bf 51}, 136--151.

\bibitem{AK1996}
Asmussen, S. and Kl\"uppelberg, C. (1996)
Large deviations results for subexponential tails,
with applications to insurance risk.
{\it Stoch.\ Proc.\ Appl.} {\bf 64}, 103--125.

\bibitem{APP}
Avram, F., Palmowski, Z. and Pistorius, M. (2008)
Exit problem of a two-dimensional risk process from a cone:
exact and asymptotic results,
{\it Ann. Appl. Probab.} {\bf 18}, 2421--2449.

\bibitem{APP2}
Avram, F., Palmowski, Z. and Pistorius, M. (2008)
A two-dimensional ruin problem on the positive quadrant,
{\it Insurance Math. Econom.} {\bf 42}, 227--234.

\bibitem{badila}
Badila, S., Boxma, O., Resing, J. and Winands, E.M.M. (2014)
Queues and risk models with simultaneous arrivals,
{\it Adv. Appl. Probab.} {\bf 46}, 812--831.

\bibitem{CWW2013_1}
Chen, Y., Wang, Y. and Wang, K. (2013)
Asymptotic results for ruin probability of a
two-dimensional renewal risk model,
{\it Stoch. Anal. Appl.} {\bf 31}, 80--91.

\bibitem{CWW2013_2}
Chen, Y., Wang, Y. and Wang, K. (2013)
Uniform asymptotics for the finite-time ruin probabilities of two kinds
of nonstandard bidimensional risk models,
{\it J. Math. Anal. Appl.} {\bf 401}, 114--129.

\bibitem{CYN2011}
Chen, Y., Yuen, K. and Ng, K. (2011)
Asymptotics for the ruin probabilities of a two-dimensional
renewal risk model with heavy-tailed claims,
{\it Appl. Stoch. Models Bus. Ind.} {\bf 27} (3), 290--300.

\bibitem{DFK} Denisov, D., Foss, S. and Korshunov, D. (2010)
Asymptotics of randomly stopped sums in the presence of heavy tails,
{\it Bernoulli} {\bf 16(4)}, 971--994.
%Dembo, A. and Zeitouni, O. (1998) {\it Large Deviations Techniques and Applications,}
%Springer-Verlag, New-York.

%\bibitem{EmbrVeraver}
%Embrechts, P. and Veraverbeke, N. (1982) Estimates for the
%probability of ruin with special emphasis on the probability of
%large claims. {\it Insurance Math. Econom.} {\bf 1}, 55--72.

\bibitem{Embrechts}
Embrechts, P., Kl\"uppelberg, C., and Mikosch, T. (1997)
{\it Modelling Extremal Events for Insurance and Finance,} Springer, Berlin.

\bibitem{FZ}
Foss, S. and Zachary, S. (2003)
The maximum of a random time interval of a random walk with long-tailed increments and negative drift,
{\it Ann. Appl. Probab.} {\bf 13}, 37--53.

\bibitem{FPZ}
Foss, S., Palmowski, Z. and Zachary, S. (2005)
The probability of exceeding a high boundary on a random time
interval for a heavy-tailed random walk,
{\it Ann. Appl. Probab.} {\bf 3}, 1936--1957.

\bibitem{FKZ}
Foss, S., Korshunov, D. and Zachary, S. (2013)
{\it An Introduction to Heavy-tailed and Subexponential Distributions,}
Springer.

\bibitem{hu}
Hu, Z. and Jiang, B. (2013)
On joint ruin probabilities of a two-dimensional risk model
with constant interest rate,
{\it J. Appl. Probab.} {\bf 50}, 309--322.

\bibitem{add2} 
Jiang T., Wang Y., Chen Y. and Xu H. (2015) Uniform asymptotic estimate
for finite-time ruin probabilities of a time-dependent bidimensional
renewal model, {\it Insurance Math. Econom.} {\bf 64}, 45--53.

\bibitem{KL2016}
Konstantinides, D. G. and Li, J. (2016)
Asymptotic ruin probabilities for a multidimensional renewal risk
model with multivariate regularly varying claims,
{\it Insurance Math. Econom.} {\bf 69}, 38--44.

\bibitem{Dima}
Korshunov, D. (2016)
On subexponential tails for the maxima
of negatively driven compound renewal and L\'evy processes,
{\it https://arxiv.org/abs/1608.09004v2.}

\bibitem{LLT2007}
Li, J., Liu, Z. and Tang, Q. (2007)
On the ruin probability of a bidimensional perturbed risk model,
{\it Insurance Math. Econom.} {\bf 41}, 185--195.

\bibitem{Mandjes}
Lieshout, P.M.D. and Mandjes, M. (2007) Tandem Brownian queues,
{\it Mathematical Methods of Operations Research} {\bf 66}, 275--298.

\bibitem{LZ2016}
Lu, D. and Zhang, B. (2016)
Some asymptotic results of the ruin probabilities in a
two-dimensional renewal risk model with some strongly
subexponential claims,
{\it Statist. Probab. Lett.} {\it 114}, 20--29.

%\bibitem{PalPist}
%Palmowski, Z. and Pistorius, M. (2008)
%The probability of exceeding a piecewise deterministic barrier
%by the heavy-tailed renewal compound process,
%{\it http://arxiv.org/abs/0805.1631}.

\bibitem{Rolski}
Rolski, T., Schmidli, H., Schmidt, V. and Teugles, J.L. (1999)
{\it Stochastic processes for insurance and finance,}
John Wiley and Sons, Inc., New York.

\bibitem{YGW2006}
Yuen, K., Guo, J., Wu, X. (2006)
On the first time of ruin in the bivariate compound Poisson model,
{\it Insurance Math. Econom.} {\bf 38}, 298--308.

\bibitem{add1}
Wang, Y.,  Cui, Z., Wang, K. and Ma, X. (2012)
Uniform asymptotics of the finite-time ruin probability for all times,
{\it J. Math. Anal. Appl.} {\bf 390}, 208--223.

\end{thebibliography}
\end{document}